\documentclass{amsart}
 \usepackage{graphicx}
\usepackage{amssymb}
\usepackage{amsmath}

\newtheorem{tw}{Theorem}

\newcommand{\cal}[1]{\mathcal{#1}}
\newcommand{\bez}{\setminus}
\newcommand{\sig}{\sigma}
\newcommand{\eps}{\varepsilon}
\newcommand{\fal}[1]{\widetilde{#1}}
\newcommand{\gen}[1]{\langle #1 \rangle}
\newcommand{\map}[3]{#1\colon #2\to #3}
\newcommand{\Map}[2]{#1\colon #2\to #2}
\newcommand{\field}[1]{\mathbb{#1}}
\newcommand{\zz}{\field{Z}}
\newcommand{\rr}{\field{R}}
\newcommand{\lst}[2]{{#1}_1,\dotsc,{#1}_{#2}}

\newcommand{\Mpm}{{\cal{M}}_{g}^{\pm}}
\newcommand{\Mg}{{\cal{M}}_{g}}
\newcommand{\Mgh}{{\cal{M}}_{0,2g+2}}
\newcommand{\Mghpm}{{\cal{M}}_{0,2g+2}^{\pm}}
\newcommand{\Mtw}{{\cal{M}}_{2}}

\newcommand{\Mtwpm}{{\cal{M}}_{2}^{\pm}}
\newcommand{\Mone}{{\cal{M}}_{1}}
\newcommand{\Mh}{{\cal{M}}^{h}_{g}}
\newcommand{\Mhpm}{{\cal{M}}^{h\pm}_{g}}

\DeclareMathOperator{\Aut}{Aut}%
\DeclareMathOperator{\Out}{Out}%
\title{The extended mapping class group is generated by three symmetries}
\author{Micha\l\ Stukow}
\thanks{Supported by BW 5100-5-0080-3}
\email{trojkat@math.univ.gda.pl}
\address{Institute of Mathematics, University of Gda\'nsk, Wita Stwosza 57,
80-952 Gda\'nsk, Poland }

\begin{document}
\begin{abstract}
We prove that for $g\geq 1$ the extended mapping class group is generated by three orientation
reversing involutions.
\end{abstract}

\maketitle



\section{Introduction}
Let $S_g$ be a closed orientable surface of genus $g$. Denote by $\Mpm$ the \emph{extended mapping
class group} i.e. the group of isotopy classes of homeomorphisms of $S_g$. By $\Mg$ we denote the
\emph{mapping class group} i.e. the subgroup of $\Mpm$ consisting of orientation preserving maps.
We will make no distinction between a map and its isotopy class, so in particular by the order of a
homeomorphism $\Map{h}{S_g}$ we mean the order of its class in $\Mpm$.

By $C_i,U_i,Z_i$ we denote the right Dehn twists along the curves $c_i,u_i,z_i$ indicated in Figure
\ref{r1}. It is known
that this set of
generators of $\Mg$ is not minimal, and a great deal of attention has been paid to the problem of
finding a minimal (or at least small) set of generators or a set of generators with some additional
property. For different approaches to this problem see
\cite{Bre-Farb,Hump,Kork1,MacMod,McCarPap1,Wajn} and references there. The main purpose of this
note is to prove that for $g\geq 1$ the extended mapping class group $\Mpm$ is generated by three
symmetries, i.e. orientation reversing involutions. This generalises a well known fact for ${\cal
M}_1^\pm\cong \text{GL}(2,\zz)$.


As was observed in \cite{GG-MS}, the fact that $\Mpm$ is generated by symmetries is rather simple.
Namely, suppose that $S_g$ is embedded in $\rr^3$ as shown in Figure \ref{r1}. Define the
\emph{sandwich symmetry} $\Map{\tau}{S_g}$ as a reflection across the $yz$-plane. Now if $u$ is any
of the curves indicated in Figure \ref{r1}, then the twist $U$ along this curve satisfies the
relation: $\tau U\tau=U^{-1}$, i.e. the element $\tau U$ is a symmetry. This proves that each of
generating twists is a product of two symmetries. Note that for the composition of mappings we use
the following convention: $fg$ means that $g$ is applied first.


\section{Preliminaries}
\begin{figure}[h]
\begin{center}
\includegraphics[width=0.7\textwidth]{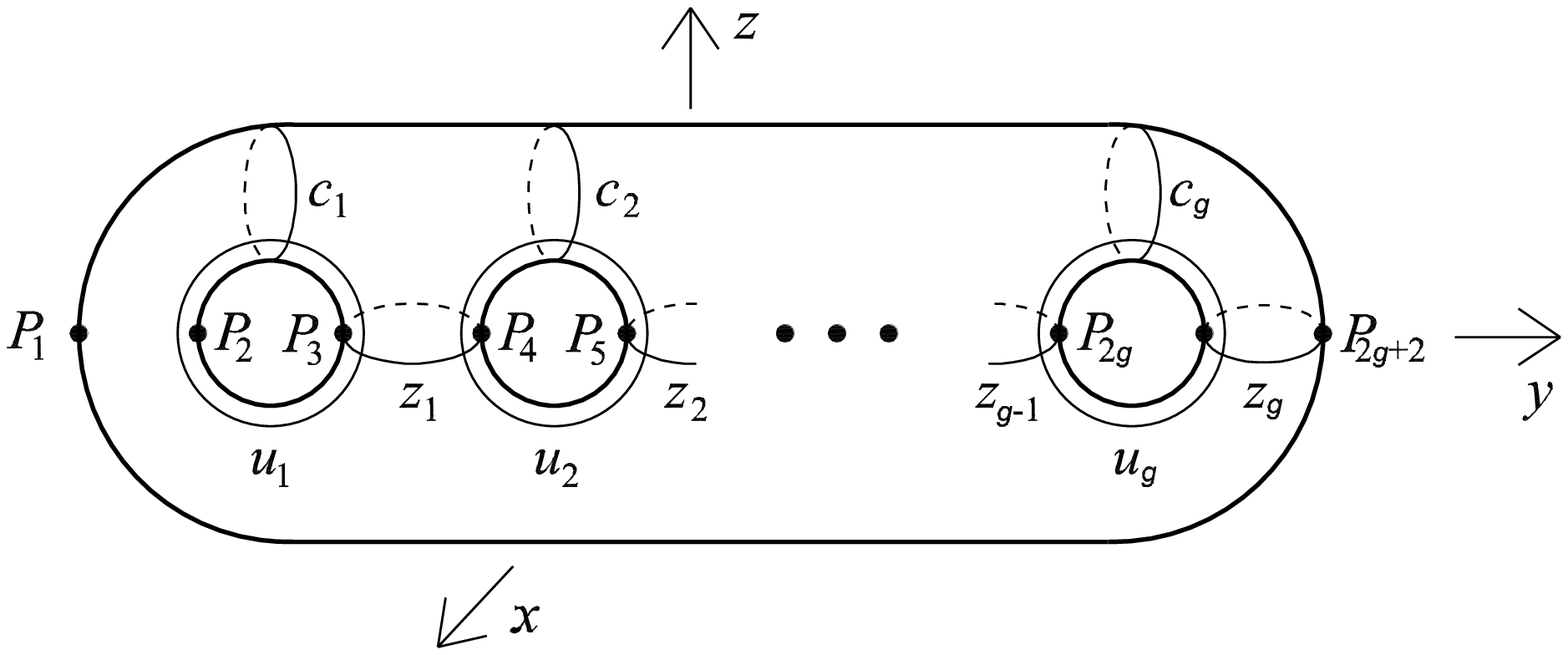}
\caption{Surface $S_g$ embedded in $\rr^3$.}\label{r1} %
\end{center}
\end{figure}

Suppose that $S_g$, for $g\geq 2$, is embedded in $\rr^3$ as shown in Figure \ref{r1}. Let
$\Map{\rho}{S_g}$ be a \emph{hyperelliptic involution}, i.e. the half turn about $y$-axis.

The \emph{hyperelliptic mapping class group} $\Mh$ is defined to be the centraliser of $\rho$ in
$\Mg$. By \cite{Bir-Hil} the quotient $\Mh/\gen{\rho}$ is isomorphic to the mapping class group
$\Mgh$ of a sphere $S_{0,2g+2}$ with $2g+2$ marked points $\lst{P}{2g+2}$. This set of marked
points corresponds (under the canonical projection) to fixed points of $\rho$ (Figure \ref{r1}). In
a similar way, we define the \emph{extended hyperelliptic mapping class group} $\Mhpm$ which
projects onto the extended mapping class group $\Mghpm$ of $S_{0,2g+2}$. Denote this projection by
$\map{\pi}{\Mhpm}{\Mghpm}$. In case $g=2$ it is known that $\Mtw={{\cal{M}}^{h}_2}$ and
$\Mtwpm={{\cal{M}}^{h\pm}_2}$.

Denote by $\sig_1,\sig_2,\ldots,\sig_{2g+1}$ the images under $\pi$ of twist generators
$C_1,U_1,Z_1,U_2,Z_2,\ldots,U_g,Z_g$ respectively. These generators of $\Mgh$ are closely related
to Artin braids, cf. \cite{Bir-Hil}.

\begin{figure}[h]
\begin{center}
\includegraphics[width=0.3\textwidth]{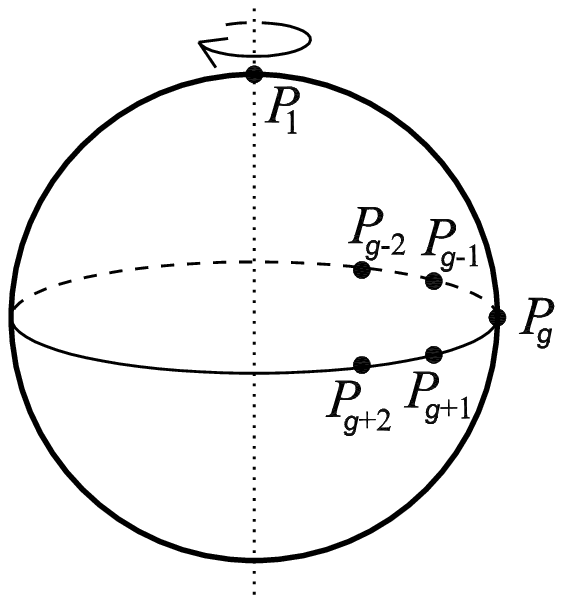}
\caption{Rotation $\fal{M}$.}\label{r2}
\end{center}
\end{figure}

Let $\Map{\fal{M}}{S_{0,2g+2}}$ be a rotation of order $2g+1$ with a fixed point $P_1$ such that:
$\fal{M}(P_i)=P_{i+1}$, for $i=2,\ldots 2g+1$ and $\fal{M}(P_{2g+2})=P_2$ (Figure \ref{r2}).
In terms of the generators $\lst{\sig}{2g+1}$ we have:
\begin{equation}
\fal{M}=\sig_2\sig_3\cdots\sig_{2g+1}.\label{eq:2}
\end{equation}
If $M'\in\Mg$ is the lifting of $\fal{M}$ of order $2g+1$, then $M=\rho M'$ is the lifting of $\fal
M$ for which $M^{2g+1}=\rho$. In particular $M$ has order $4g+2$. Using the technique described in
\cite{McCarPap1} it is easy to write $M$ as a product of twists: $M=U_1Z_1U_2Z_2\cdots U_gZ_g.$

Since every finite subgroup of $\Mg$ can be realised as the group of automorphisms of a Riemann
surface \cite{Kerk2}, $M$ has maximal order among torsion elements of $\Mg$ \cite{Wiman2}.
Geometric properties of $M$ played a crucial role in the problem of finding particular sets of
generators for $\Mg$ and $\Mpm$ cf. \cite{Bre-Farb,Kork1,MacMod,Wajn}.

\begin{figure}[h]
\begin{center}
\includegraphics[width=0.6\textwidth]{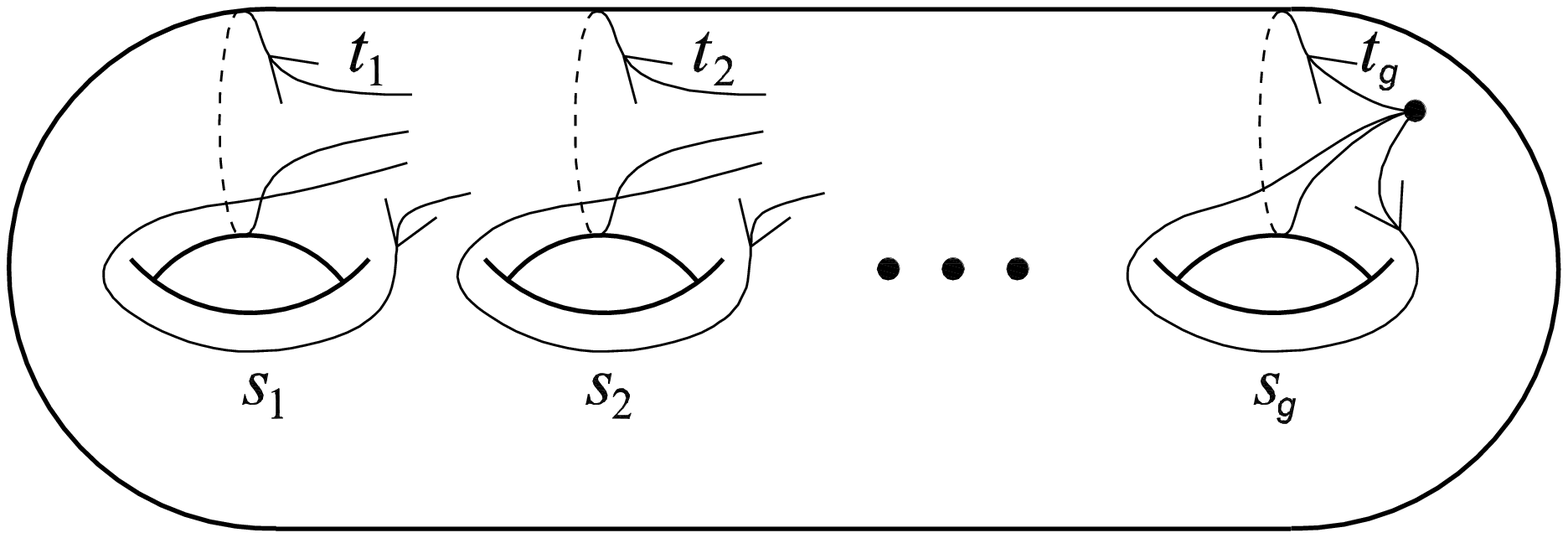}
\caption{Generators of $\pi_1(S_g)$.}\label{r3} %
\end{center}
\end{figure}

Following \cite{Bir2}, let $t_1,s_1,\ldots,t_g,s_g$ be generators of the fundamental group
$\pi_1(S_g)$ as in Figure \ref{r3}.
In terms of these generators, $\pi_1(S_g)$ has the single defining relation:
$R=s_g^{t_g}s_{g-1}^{t_{g-1}}\cdots s_1^{t_1}s_1^{-1}s_2^{-1}\cdots s_g^{-1}$, where by $a^b$ we
denote the conjugation $bab^{-1}$.

 It is well known \cite{MKS} that the mapping class group $\Mpm$ is isomorphic to the group $\Out(\pi_1(S_g))$ of
 outer automorphisms of $\pi_1(S_g)$. In terms of this isomorphism, elements of $\Mg$
 correspond to the elements of $\Out(\pi_1(S_g))$ which map the relation $R$ to its conjugate,
 and elements of $\Mpm\bez\Mg$ to those elements of $\Out(\pi_1(S_g))$ which map $R$ to a
 conjugate of $R^{-1}$.

 Using representations of twist generators as automorphisms of $\pi_1(S_g)$ \cite{Bir2} we could
 derive the following representation for the rotation $M$:
 \begin{alignat*}{3}
M:\ &t_i &&\mapsto s_i^{t_i}\cdots s_1^{t_1}t_1& &\quad\text{for } i=1,\ldots, g\\
&s_i &&\mapsto t_1^{-1}s_1^{-t_1}\cdots s_i^{-t_i}t_{i+1}t_i^{-1} s_i^{t_i}\cdots s_1^{t_1}t_1
& &\quad\text{for } i=1,\ldots, g-1 \\
&s_g &&\mapsto t_1^{-1}s_1^{-t_1}\cdots s_g^{-t_g}t_g^{-1} s_g^{t_g}\cdots s_1^{t_1}t_1.& &
\end{alignat*}
As in the case of maps and their isotopy classes, we abuse terminology by identifying an element of
$\Out(\pi_1(S_g))$ with its representative in $\Aut(\pi_1(S_g))$.

\section{$\Mpm$ is generated by 3 symmetries}
If we represent the action of the rotation $\fal{M}$ as the orthogonal action on the unit sphere,
it becomes obvious that $\fal{M}$ can be written as a product of two symmetries. To be more
precise, if $\fal{\eps}_1$ is the symmetry across the plane passing through $P_1,P_{g}$ and the
center of the sphere (Figure \ref{r2}), then $\fal{M}=\fal{\eps}_1\fal{\eps}_2$, where
$\fal{\eps}_2$ is another symmetry.

Tedious but straightforward computations show that one of the liftings $\eps_1\in \Mpm$ of
$\fal{\eps}_1$ has the following representation as an automorphism of $\pi_1(S_g)$:
 \begin{alignat*}{1}
\eps_1:\ &t_i \mapsto\ t_{g-1}^{-1}s_1^{-1}\cdots s_{g-1-i}^{-1},\quad s_i \mapsto\ t_{g-1-i}^{-1}t_{g-i}\quad\text{for } i=1,\ldots, g-2\\
&t_{g-1} \mapsto\ t_{g-1}^{-1},\quad s_{g-1} \mapsto\ s_{g}\cdots s_1 t_1,\quad t_{g} \mapsto\
t_{g-1}^{-1}t_g,\quad s_{g} \mapsto\ s_g^{-1}
\end{alignat*}
To obtain the above representation we proceed as follows: take a generator $u$ of $\pi_1(S_g)$,
find the image $\fal{u}$ of $u$ under projection $S_g\rightarrow S_{0,2g+2}$, find
$\fal{\eps}_1(\fal{u})$, lift back $\fal{\eps}_1(\fal{u})$ to $S_g$ and finally express the
obtained loop as a product of generators $t_1,s_1,\ldots,t_g,s_g$ of $\pi_1(S_g)$.

We would like to point out that although the above procedure is a bit subtle, it is quite simple to
verify that the obtained formulas are correct. In fact, it is enough to check that $\eps_1^2=1$ and
$\eps_1(R)$ is conjugate to $R^{-1}$. Moreover, the representation of $\eps_2=\eps_1M$ is given by
the
 following formulas:
\begin{alignat*}{1}
\eps_2:\ &t_i \mapsto\ (t_{g-1}^{-1}s_1^{-1}\cdots
s_{g-1-i}^{-1}t_{g-1-i}^{-1})(s_{g-i}^{-t_{g-i}}\cdots
s_{g-1}^{-t_{g-1}})t_{g-1}s_{g-1} \quad\text{for } i=1,\ldots, g-2\\
&t_{g-1} \mapsto\ t_{g-1}^{-1}s_g^{t_g}t_{g-1}s_{g-1},\quad t_{g} \mapsto\ s_{g-1}\\
&s_{i} \mapsto\ s_{g-1}^{-1}t_{g-1}^{-1}(s_{g-1}^{t_{g-1}}\cdots s_{g-i}^{t_{g-i}})
(s_{g-1-i}^{t_{g-1-i}})(s_{g-i}^{-t_{g-i}}\cdots s_{g-1}^{-t_{g-1}})t_{g-1}s_{g-1} \quad\text{for }
i=1,\ldots, g-2\\
&s_{g-1} \mapsto\ (s_{g-1}^{-1}{t_{g-1}^{-1}}s_{g}^{-t_{g}})t_g(s_g^{t_g}t_{g-1}s_{g-1}),\quad
s_{g} \mapsto\ s_{g-1}^{-1}t_g^{-1}t_{g-1}s_{g-1}
\end{alignat*}
It is straightforward to verify that $\eps_2^2$ is an identity in $\Out(\pi_1(S_g))$.
\begin{tw}
For each $g\geq 1$, the extended mapping class group $\Mpm$ is generated by three symmetries.
\end{tw}
\begin{proof}
As observed in the introduction, the result is well known for $g=1$, but for the sake of
completeness let us prove this in more geometric way. Since $\Mone=\gen{U_1,C_1}$ (Figure \ref{r1})
and $\tau U_1\tau=U_1^{-1}$, $\tau C_1\tau=C_1^{-1}$, the group ${\cal M}_1^{\pm}$ is generated by
the symmetries $\tau,\tau U_1,\tau C_1$.

Now suppose that $g\geq 2$. Let $\eps_1$ and $\eps_2=\eps_1M$ be the symmetries defined above.
Since $\eps_1(t_{g-1})=t_{g-1}^{-1}$ we have $\eps_1C_{g-1}\eps_1=C_{g-1}^{-1}$, i.e.
$\eps_3=\eps_1C_{g-1}$ is a symmetry. In particular
$\gen{\eps_1,\eps_2,\eps_3}\supset\gen{\eps_1\eps_2,\eps_1\eps_3}=\gen{M,C_{g-1}}$. But by
\cite{Kork1} the latter group is equal to $\Mg$. Since $\gen{\eps_1,\eps_2,\eps_3}$ contains
orientation reversing element, this proves that $\gen{\eps_1,\eps_2,\eps_3}=\Mpm$.
\end{proof}




\section*{Acknowledgements}
The author wishes to thank the referee for his/her helpful suggestions.

\end{document}